\documentclass[11pt]{amsart}
\usepackage{amsmath,amsthm, amssymb, amscd,amsfonts}
\usepackage{mathrsfs}
\usepackage{enumerate}
\makeatletter
\@namedef{subjclassname@2010}{%
  \textup{2010} Mathematics Subject Classification}
\makeatother
\newtheorem{lemma}{Lemma}[section]
\newtheorem{theorem}{Theorem}[section]

\newtheorem{corollary}{Corollary}[section]
\newtheorem{definition}{Definition}[section]
\newtheorem{remark}{Remark}[section]

\begin{document}

\title{CONTINUOUS $k$-FRAMES AND THEIR DUALS} 
 
\author[ Gh.Rahimlou] {Gholamreza Rahimlou}
\address[First Author]{Department of Mathematics, Shabestar Branch,Islamic Azad University, Shabestar,Iran.}
\email{grahimlou@gmail.com}

\author[R. Ahmadi]{Reza Ahmadi}
\address[Second Author]{Institute of Fundamental Science, University of Tabriz, Tabriz, Iran.}
\email{rahmadi@tabrizu.ac.ir}

\author[M. A. Jafarizadeh]{Mohammad Ali Jafarizadeh}
\address[Third Author]{Faculty of Physic,  University of Tabriz, Tabriz, Iran.}
\email{jafarizadeh@tabrizu.ac.ir}

\author[S.Nami]{Susan Nami}
\address[Fourth Author]{Faculty of Physic,  University of Tabriz, Tabriz, Iran. }
\email{ S.Nami@tabrizu.ac.ir}

 \begin{abstract}
$k$-frames were recently introduced by G\u avru\c ta in Hilbert  spaces to study atomic systems with respect to a bounded linear operator. A continuous frame is a family of vectors in Hilbert space which allows reproductions of arbitrary elements by continuous super positions. We construct a continuous $k$-frame, so called c$k$-frame along with an atomic system for this version of frames. Also we introduce a new method for obtaining the dual of a c$k$-frame and we prove some new results about it. 
\end{abstract}



 \maketitle
 
\section{Introduction}

Frames were first introduced in the context of non-harmonic Fourier
series \cite{hy44}. Outside of signal processing, frames did not
seem to generate much interest until the ground breaking work of
\cite{ga94}. Since then the theory of frames began to be more widely
studied. During the last $20$ years the theory of frames has been
growing rapidly and several new applications have been developed. For
example, besides traditional application as signal processing, image
processing, data compression, and sampling theory, frames are now
used to mitigate the effect of losses in pocket-based communication
systems and hence to improve the robustness of data transmission 
\cite{dt97}, and to design high-rate constellation with full
diversity in multiple-antenna code design \cite{hir98}. In
\cite{aklr, am00, bbk} some more applications have been developed.

 The structure of this article is as follows; in Section 1 we review some basic properties of frame theory in Hilbert spaces. In Section 2 we introduce c$k$-frames and some fundamental properties about them are discussed. Finally, in Section 3 we introduce a new method for obtaining the dual of a c$k$-frame and we prove some new results about it.

Throughout the paper, $H$ and $H_0$ are Hilbert spaces, ${(H_0)_{1}}$ is the closed unit ball in $ H_0$ , $(X,\mu)$ is a $\sigma$-finite measure space, $\mathcal{L}(H_0,H)$ is the set of all linear mappings of $H_0$ to $H$ and  $\mathcal{B}(H_0, H)$ is the set of all bounded linear mappings. Instead of $\mathcal{B}(H,H)$, we simply write $\mathcal{B}(H)$. 
 \begin{definition}\label{definition1}
 Let $f:X\rightarrow H$ be a weakly
measurable (i.e., for all $h\in H$, the mapping $x
\mapsto\langle f(x),h\rangle$ is measurable). Then $f$ is called a 
c-frame  for $H$ if there exist $0 \leq A\leq B <\infty$ such that for
all $h\in H$,
\begin{eqnarray*}
 A\|h\|^2\leq  \int_{X}|\langle f(x),h\rangle|^2d\mu \leq B\|h\|^2.
 \end{eqnarray*}
\end{definition}
The constants $A$ and $B$ are called c-frame bounds. If $A,B$ can be
chosen so that $A=B$, we call this c-frame an $A$-tight frame,  and if
$A=B=1$ it is called a c-Parseval frame. If we only have the upper
bound, we call $f$ a c-Bessel mapping for $H$.
The representation space employed in this setting is
$$L^2(X,H)=\{\varphi:X\rightarrow H | \ \ \varphi  \mbox{ is measurable and} \ \  \|\varphi\|_2 <\infty\},$$ where
$\|\varphi\|_2=(\int_X||\varphi(x)||^2d\mu)^{\frac{1}{2}}$. 

 For each $ f,g\in L^{2}(X,H) ,$ the mapping $ x\rightarrow\langle f(x),g(x)\rangle$
of $X$ to   $\mathbb{C}$  is measurable, and it be can proved that $L^{2}(X,H) $ is a Hilbert spase with the inner product defined by $$\langle f,g\rangle_{L^{2}}=\int_{X}\langle f(x),g(x)\rangle\,d\mu. $$

We shall write $L^2(X)$ when $H= \mathbb{C}$.

\begin{theorem}[\cite{hy41}] \label{t11}
Let  $f:X \to H$ be a c-Bessel mapping  for $H$, and $u\in \mathcal{B}(H , H_0)$. Then $uf : X \to  H_0$ is a c-Bessel mapping for $H_0$ with $$ uT_f=T_{uf}$$.
\end{theorem}

\begin{theorem}[\cite{do}]\label{x1}
Suppose $ H,H_{1},H_{2} $ are Hilbert spaces, $L_{1}\in \mathcal{B}(H_{1},H)$ and $ L_{2}\in \mathcal{B}(H_{2},H)$. Then the following  assertions are equivalent:
\item
$\mathcal{R}(L_{1})\subset \mathcal{R}(L_{2})$ ,
\item
$\exists \lambda \geq 0 ,\hspace {2mm} \text{such that} \hspace {2mm} L_{1}L_{1}^{*}\leq \lambda L_{2}L_{2}^{*}$,
\item
$\text{Tere exist} \hspace {1mm} X\in \mathcal{ B}(H_{1},H_{2}) \hspace {1mm}\text{such that} \hspace {2mm} L_{1}=L_{2}X$.
\end{theorem}

\section{c$k$-frames}
\begin{definition}
  Let $ H_{0}\subseteq H$ and $k\in \mathcal{B}(H_0, H)$.  A Bessel sequence $ \lbrace f_{n}\rbrace \subseteq H $ is called a family of local $k$-atoms for $ H_{0} $ if there exist a sequence $ \lbrace c_{n}\rbrace $ of linear functionals on $ H_{0} $ such that the following conditions are satisfied:
\item
 There exists $ c>0 $ such that for each $ f\in H_{0}$
\begin{equation*}
\sum_{n}| c_{n}(f)|^{2}\leq c||f||^{2}, and
\end{equation*}
\item
$$kf=\sum_{n} c_{n}(f) f_{n}, \hspace{2mm} f\in H_{0}.
$$
 In this case we say that the pair $ \{f_{n},c_{n}\}$ provides an $k$-atomic decomposition for $ H_{0} $.\\
If $k$ is the identity map, then we say that $ \lbrace f_{n}\rbrace \subseteq H $ is a family of local atoms for $H_{0}.$
\end{definition}  

\begin{definition}
 Let  $k\in \mathcal{B}(H_0,H)$,  and $ \lbrace f_{n}\rbrace \subseteq H $.  We say that  the sequence $ \lbrace f_{n}\rbrace $ is a k-frame for $ H $ with respect to $H_0$,  if there exist  constants $ A,B>0 $ such that 
\begin{equation*}
A||k^{*}h||\leq \sum_{n}|\langle h, f_{n}\rangle |^{2}\leq B ||h||^{2}, \ \  h \in H.
\end{equation*}
\end{definition}
\begin{definition}
Let $ f:X\rightarrow H$ be weakly measurable. We define the map $ \int_{X} \cdot fd\mu: L^{2}(X)\rightarrow H $ as follow: 
\begin{equation*}
\langle \int_{X}gfd\mu,h \rangle:=\int_{X}g(x)\langle f(x),h \rangle d\mu ,\hspace {2mm} h\in H.
\end{equation*}

It is clear that, the vector valued integral $ \int_{X}gfd\mu $ exists in $H$ if for each $ h\in H,$ $ \int_{X}g(x)\langle f(x),h\rangle d\mu$ exists. 
\end{definition}

\begin{lemma}[\cite{ir97}]\label{l1}
Let $f:X\rightarrow H$ be weakly measurable. For each $ g\in L^{2}(X),$ the value of $\int_{X}gfd\mu $ exists in $H$ if and only if  for each $ h \in H$, $ \langle f,h \rangle \in L^{2}(X) $.
\end{lemma}

\begin{lemma}\label{x4}
Let $f:X\rightarrow H $ be weakly measurable.  Then $ f $ is a c-Bessel mapping for $H$ if and only if  for each $ g\in L^{2}(X) $,  $ \int_{X}gfd\mu $ exists in $H$.
\end{lemma}
\begin{proof}
Suppose that  for each $ g\in L^{2}(X) $,  $ \int_{X}gfd\mu $ exists.  By Lemma \ref{l1}, for each $ h\in H$, $ \langle f,h \rangle \in L^{2}(X) $.  We have
\begin{align*}
 | \int_{X}gfd\mu | &= \sup_{||t||=1}| \langle \int_{X}gfd\mu , t\rangle  |\\
 &= \sup_{||t||=1}  | \int_{X}g(x)\langle f(x),t\rangle d\mu  | \\
 &\leq  || g  ||_{2} \sup_{||t||=1}   ||\langle f,t\rangle   ||_{2}.
\end{align*}
Since for every $ x\in X $,
\begin{equation*}
\sup_{||t||=1}  | \langle f(x),t\rangle  |\leq  | |f(x) | | < \infty ,
\end{equation*}
so  by Banach-Steinhaus theorem,
$
\sup_{||t||=1}  | |\langle f,t\rangle | |<\infty. $
Hence 
$$||\int_X \cdot f d\mu || \leq \sup_{||t||=1} ||\langle f, t  \rangle ||_2 < \infty.$$
 The above inequality implies that $ \int_{X} \cdot fd\mu $ is bounded and  $\displaystyle \sup_{|| t||=1} || \langle f,t\rangle ||_{2} $ is an upper bounded for $ \int_{X} \cdot fd\mu. $ 
Now the adjoint of $ \int_{X}.fd\mu $ is calculated as follow:\\  For each $h\in H$ and $g\in L^{2}(X)$,  we have
 \begin{align*} \langle g, (\int_X \cdot f d\mu )^*(h) \rangle &= \langle\int_Xg fd\mu, h \rangle \\&=\int_Xg(x) \langle f(x), h \rangle d\mu  \\&=\langle g,\langle h, f\rangle \rangle.
\end{align*}
Thus, for each  $ h\in H,$
 \begin{equation}\label{e1}
( \int _{{X}} \cdot fd\mu ) ^{*}(h)=\langle h,f\rangle. \\    
 \end{equation}\\ \\ 
Therefore,
\begin{align*} \int_X |\langle h,f(x) \rangle|^2 d\mu &=||(\int_X \cdot f d\mu )^*(h)||^2 \\&\leq ||(\int_X \cdot f d\mu)^*||^2 ||h||^2\\& =||\int_X \cdot f d\mu||^2 ||h||^2  \\& \leq (\sup_{||t||=1} ||\langle f, t \rangle||_2^2)  ||h||^2.\end{align*}
Hence, $ f $ is a c-Bessel mapping for $H$. Now, if $ f $ is a c-Bessel mapping for $H$  then for each $h\in H$ we have $\langle h,f\rangle \in L^{2}(X) $.  Consquently by lemma \ref{l1} for each $ g  \in L^{2}(X) $,  $ \int_{X}.gfd\mu $ exists.
\end{proof}

\begin{remark}
Let $f:X\rightarrow H$ be a c-Bessel mapping for $H$. The synthesis operator  is
defined by
$$
T_f:L^2(X)\rightarrow H ,  \
 T_f(g)=\int_X g f\,d\mu.
$$
 Hence, for each  $g \in L^2(X)$,  and $h\in H$,
$$\langle\int_{X}gf\,d\mu,h\rangle=\int_{X}g(x)\langle f(x),h\rangle\,d\mu.$$
By \eqref{e1} the analysis operator is defined by
$$
T^{\ast}_f:H\rightarrow L^{2}(X), \ \ 
T^{\ast}_{f}(h)=\langle h,f\rangle.$$
So, for the frame operator $S_f:=T_fT^*_f$ we have 
\begin{align*}
S_f(h)=\int_{X}\langle h,f\rangle f\,d\mu, \ \ h\in H,
\end{align*}
and
\begin{align*}
\langle S_f(h),h\rangle=\Vert T^{\ast}_f(h)\Vert^2
=\int_{X}\vert\langle h,f\rangle\vert^2\,d\mu,\ \  h\in H.
\end{align*}
\end{remark}
\begin{definition} 
Let   $  H_{0}\subseteq H   $ . Suppose that  $f: X\rightarrow  H $ is  weakly measurable and $  k \in \mathcal{B}(H_{0},H) $.  Then f is called a family of  local  c$k$-atoms for    $  H_{0}$  if the following conditions are satisfied:
For each $ g\in L^{2}(X),$ the vector valued integral $~\int_{X}gfd\mu $ exists in $H$.
There exist some $ a >0 $ and $ \ell:X\rightarrow \mathcal{L}(H_{0},)$
such that for each $ h\in H_{0},~\ell (\cdot) (h)\in L^{2}(X) $ and  also
\begin{align*}
|| &\ell (\cdot)(h)||_{2}\leq a|| h || ,\\ & k h=\int_{X}\ell (\cdot)(h)fd\mu .
\end{align*} 
If $ k $  is the identity function on $ H_{0} $ then $ f $ is called a family of local  atoms  for $H_{0} $.
\end{definition}

\begin{definition}
Let $ k \in \mathcal{B}(H_{0},H) $ and $ f:X\rightarrow H $ be weakly measurable.  Then  the map $ f $ is called a ck-frame with respect to  $ H_{0} $, if there exist constants  $ A,B>0 $ such that for each $h\in H,$
\begin{equation*}
A|| k^{*} h||^{2}\leq \int_{X}|\langle h,f(x)\rangle |^{2}d\mu \leq B|| h||^{2}.
\end{equation*}
\end{definition}
\begin{theorem}\label{t5}
 Let $ H_0\subseteq H $. Let $ f:X\rightarrow H $ be weakly measurable, and $ k \in \mathcal{B}(H_{0},H) $. Then the following assertions are equivalent:
 \item
$ f $ is a family of local   ck-atoms for $ H_{0} $.\\
\item  f is a   ck-frame for $H$ with respect to $ H_{0} $.\\
\item $ f $ is a c-Bessel mapping for $ H $,  and there exist $ g \in \mathcal{B}(H_{0},L^{2}(X)) $ such that
\begin{equation*}
k h=\int_{X}g(h)fd\mu ,\quad h\in H_{0}.
\end{equation*} 
\end{theorem}
\begin{proof}
$(1)\Rightarrow (2)$. 
 By the hypothesis and Lemma \ref{x4}, $ f $ is a c-Bessel mapping for $H$. For each  $ h\in H $ we have
\begin{align*}
|| k^{*} h|| &=\sup_{||t||=1}| \langle k^{*}(h),t\rangle | \\&=\sup_{||t||=1}| \langle h,k(t)\rangle.
\end{align*}
Now  by (1) there exist $ c>0 $ and $ \ell: X\rightarrow  \mathcal{L}(H_{0},\mathbb{C}) $ such that for every $ h\in H_{0} $, $ \ell (\cdot)(h)\in L^{2}(X) $,  and also
\begin{align*}
&|| \ell(\cdot)(h)||_{2}\leq c || h ||,\\& \quad k h=\int_{X}\ell(\cdot)(h)fd\mu .
\end{align*}
So for each $h\in H,$
\begin{align*}||k^*(h)||^2&= \sup_{||t||=1} |\langle h,  \int_X \ell(\cdot) (t)fd\mu \rangle |^2\\&=\sup_{||t||=1} |\int_X  \ell(x) (t)\langle h,f(x) \rangle  d\mu|^2 \\ &\leq \sup_{||t||=1}||\ell(\cdot)(t)||_2^2 (\int_X|\langle h, f(x) \rangle|^2 d\mu)  \\&\leq \sup_{||t||=1}c^2 ||t||^2 \int_X |\langle h,f(x) \rangle |^2 d\mu \\&=c^2 \int_X| \langle h, f(x) \rangle |^2 d\mu.
 \end{align*} 
$(2)\Rightarrow (3).$ Since $ f $ is a  c-Bessel mapping for $ H $,  $ T_{f}: L^{2} (X) \rightarrow H$ is a bounded linear operator.
By (2), for each $ h\in H $
\begin{equation*}
A|| k^{*}(h)||^{2}\leq || T_{f}^{*}(h)||^2 .
\end{equation*}
Now for each $ h\in H $, we have
\begin{align*}
A\langle kk^{*}(h),h\rangle =A|| k^{*}(h)||^{2}\leq || T_{f}^{*}(h)||^{2}=\langle T_{f}T_{f}^{*}(h),h\rangle .
\end{align*}
Thus
\begin{equation*}
\displaystyle kk^{*}\leq \frac {1}{A} T_{f}T_{f}^{*}.
\end{equation*}
Finally, by Theorem \ref{x1},  there exist a bounded linear operator $ M: H_{0}\longrightarrow  L^{2}(X)$ such that $ k=T_{f}M $. So for each $ h\in H_{0} $ 
\begin{equation*}
kh=T_{f}(M(h))=\int_{X}M(h)fd\mu .
\end{equation*}
$(3)\Rightarrow (1).$ Since $ f $ is a  c-Bessel mapping for $ H $, thus by Lemma \ref{x4} for each $ g\in L^{2}(X) $, $ \int_{X}gfd\mu $ exists. By (3) there exist $ g \in \mathcal{B}(H_{0} ,L^{2}(X) )$ such that
\begin{equation*}
kh=\int_{X}g(h)fd\mu ,\quad h\in H_{0}.
\end{equation*}
Now we define 
\begin{equation*}
\ell : X\rightarrow L(H_{0},),\quad \ell(\cdot)(h):=g(h)(\cdot), \ \ h\in H_0,
\end{equation*}
so we have
\begin{equation*}
k(h)=\int_{X}\ell(\cdot)(h)fd\mu ,\qquad h\in H_{0}, 
\end{equation*}
also
\begin{equation*}
|| \ell(\cdot)(h)||_{2}= || g(h)(\cdot)||_{2}\leq || g||  || h|| .
\end{equation*}
This completes the proof.
\end{proof}
\begin{theorem}\label{t6}
Let $k \in B(H_{0},H)$. Suppose that $ f: X\rightarrow H $ is weakly measurable. Then f is a ck-frame for $H$ with respect to $ H_{0} $ if and only if the map
\begin{equation*}
L_{f}:L^{2}(X)\longrightarrow H,   \\ \quad L_{f}(g)=\int_{X}gfd\mu, 
\end{equation*}
  is well-defined  bounded  linear  operator with  $ \mathcal{R}(k)\subset \mathcal{R}(L_{f}). $        
\end{theorem}
\begin{proof}
First we assume that  $ f $ is  a ck-frame for $H$ with respect to $ H_{0} $. Then by definition, $ f  $ is a  c-Bessel mapping for $H$. We have
\begin{align*}
A || k^{*}(h)||^{2}&\leq \int_{X}| \langle h,f(x)\rangle|^{2}d\mu =|| T^{*}_{f}(h)||^{2},
\end{align*}
thus
$$Akk^{*}\leq T_{f}T_{f}^{*},$$
and  Theorem \ref{x1} implies that
$$\mathcal{ R}(k)\subset\mathcal{R }(T_{f}).$$
Since $L_{f} =T_f$, so $L_f$ is a bounded linear operator. Now, let 
$$ L_{f}: L^{2}(X)\longrightarrow H$$
 $$L_{f}(g)=\int_{X} gfd\mu, $$
be a  well-defined bounded  linear operator of $ L^{2}(X) $ into $ H $ with $ \mathcal{R}(k)\subset \mathcal{R}(L_{f}).$ By Lemma \ref{x4}, $f$ is a c-Bessel mapping for $H$. So it is sufficient to show that it has a lower ck-frame bound. Since $k$ and $L_{f}$ are bounded  linear operators and $\mathcal{R}(k) \subset  \mathcal{R}(L_{f})$,  by Theorem \ref{x1},  there exists $A>0 $ such that $A kk^{*}\leq L_{f}L_{f}^{*}$.  Now for each $h\in H$,
$$\langle A kk^{*}(h),h\rangle \leq \langle L_{f}L_{f}^{*}(h),h\rangle,$$
consquently by  \eqref{e1} we have
$$A ||   k^{*}(h)  ||^{2}\leq || L_{f}^{*}(h)  ||^{2}=\int_{X} |\langle h,f(x)\rangle |^{2}d\mu,$$
and the  proof is complete.
\end{proof}
Generally, in ck-frames, as $k$-frames, the frame operator is not invertible. However, we have the following:
\begin{theorem}\label{re1}
Let $k\in \mathcal{B}(H_0,H)$, and $f:X\to H$ be a ck-frame for $H$ with respect to $H_0$,  with bounds A,B. If
  $k$ be closed range then $S_f$ is invertibale on $\mathcal{R}(k)$, and 
for each $h\in \mathcal{R}(k)$
$$ B^{-1}||h||^2 \leq \langle  (S_f|_{\mathcal{R}(k)})^{-1}h,h\rangle \leq A^{-1}\Vert k^{\dagger}\Vert^2 \Vert h\Vert^2.$$
\end{theorem}
\begin{proof}
For each $h\in H$ 
$$A\Vert k^{\ast}h\Vert^2\leq\int_{X}\vert\langle h,f\rangle\vert^2\,d\mu=\langle S_f(h),h\rangle\leq B\Vert h\Vert^2$$
hence
\begin{eqnarray*}
Akk^{\ast}\leq S_f \leq BI.
\end{eqnarray*}
Since $kk^{\dagger}\vert_{\mathcal{R}(k)}=I_{\mathcal{R}(k)}$, so for each $h \in \mathcal{R}(k)$
\begin{align*}
\Vert h\Vert=\Vert I^{\ast}_{\mathcal{R}(k)}h\Vert=\Vert (k^{\dagger}\vert_{\mathcal{R}(k)})^{\ast}k^{\ast}h\Vert \leq\Vert k^{\dagger}\Vert.\Vert k^{\ast}h\Vert.
\end{align*}
Therefore,  for each $h\in \mathcal{R}(k)$,
$$A\Vert k^{\dagger}\Vert^{-2}\Vert h\Vert^2 \leq \langle  S_f(h),h\rangle \leq B||h||^2.$$
whence, $S_f$ is invertibale on $\mathcal{R}(k)$, and 
for each $h\in \mathcal{R}(k),$
$$ B^{-1}||h||^2 \leq \langle  (S_f|_{\mathcal{R}(k)})^{-1}(h),h\rangle \leq A^{-1}\Vert k^{\dagger}\Vert^2 \Vert h\Vert^2.$$
\end{proof}
\begin{corollary}
Let $k\in \mathcal{B}(H_0,H)$, and $f:X\to H$ be a ck-frame for $H$ with respect to $H_0$,  with bounds A,B. If $k$ be closed range with $\mathcal{R}(k) \subset \mathcal{R}(f)$ then $f$ is a c-frame for $\mathcal{R}(k)$ with bounds $A\Vert k^{\dagger}\Vert^{-2}$ and  $B$, respectively.
\end{corollary}
\section{c$k$-Duals}
In this section, we introduce a dual of ck-frames and prove some theorems about them.
Throughout this section, the orthogonal projection of $H$ onto a closed subspace $V\subseteq H$ is denoted by $\pi_V$.
\begin{theorem} \label{t3}
Let $ k \in \mathcal{B}(H_{0},H)$, and let  $f:X \to H$ be c-Bessel mappings for $H$, and $g:X\to H_0$ be a c-Bessel mapping for $H_0$. Then the following assertions are equivalent:
\item
 For each $h_0\in{H_0}$, $kh_0=T_{f}( \langle h_0,g\rangle)$.
\item
For each $h\in{H}$, $k^{*}h=T_{g}( \langle h,f\rangle)$.
\item For each $h\in H,h_0\in{H_0}$,
$$\langle kh_0,h\rangle = \int_{X}\langle h_0,g(x)\rangle \langle f(x),h\rangle \,d\mu.$$
\item
For each $h\in H,h_0\in{H_0}$,
\begin{align*}
\langle k^{*}h,h_0\rangle = \int_{X}\langle h,f(x)\rangle \langle g(x),h_0\rangle \,d\mu .
\end{align*}

\item
For each orthonormal bases  $\lbrace \gamma_{j}\rbrace_{j\in{J}}$ for $H_0$, and $\lbrace e_{i}\rbrace_{i\in{I}}$  for $H,$
\begin{align*}
\langle k^{*}e_{i},\gamma_{j}\rangle = \int_{X}\langle e_{i},f(x)\rangle \langle g(x),\gamma_{j}\rangle \,d\mu \ , \ i\in{I} \ j\in{J}.
\end{align*} 
\end{theorem}
\begin{proof}
$(1)\Rightarrow (2)$. \ If $h\in H$, and  $h_0 \in H_0$ then
\begin{align*}
\langle h_0, k^{*}h \rangle&=\langle T_f(\langle h_0,g\rangle),h\rangle\\
&=\int_{X}\langle h_0,g(x)\rangle \langle f(x),h\rangle\,d\mu\\
&=\overline{\int_{X}\langle g(x),h_0\rangle \langle h,f(x)\rangle\,d\mu}\\
&=\overline{\langle T_g(\langle h,f\rangle),h_0\rangle}\\
&=\langle h_0,T_g(\langle h,f\rangle)\rangle.
\end{align*} 
Hence  (2) is proved.

$(2)\Rightarrow (1)$.  \ By the same proof of $ (1) \Rightarrow (2)$. 

$(4) \Rightarrow(5) $, $(2)\Leftrightarrow(3)$, and  
$(3)\Leftrightarrow(4)$ are evident.

$(5)\Rightarrow(4)$. \ Let $h\in H , h_0\in H_0$. Then
\begin{align*}
 \int_{X}\langle h,f(x)\rangle \langle g(x),h_{0}\rangle \,d\mu &=\Big\langle\langle h,f\rangle ,\langle h_{0},g\rangle\Big\rangle_{L^{2}}\\
&=\Big\langle\langle h,\sum_{i}\overline{\langle e_{i},f\rangle}\ e_{i}\rangle ,\langle h_{0},\sum_{j}\overline{\langle \gamma_{j},g\rangle}\ \gamma_{j}\rangle\Big\rangle_{L^{2}}\\
&=\sum_{i,j}\big\langle\langle h,\overline{\langle e_{i},f\rangle} e_{i}\rangle\langle h_{0},\overline{\langle \gamma_{j},g\rangle} \gamma_{j}\rangle \big\rangle_{L^{2}}\\
&=\sum_{i,j}\langle h,e_{i}\rangle\langle \gamma_{j},h_{0}\rangle\big\langle \langle e_{i},f\rangle, \langle \gamma_{j},g\rangle\big\rangle_{L_{2}}\\
&=\sum_{i,j}\langle h,e_{i}\rangle\langle \gamma_{j},h_{0}\rangle\int_{X}\langle e_{i},f(x)\rangle\langle g(x),\gamma_{j}\rangle\,d\mu\\
&=\sum_{i,j}\langle h,e_{i}\rangle\langle \gamma_{j},h_{0}\rangle\langle k^{*}e_{i},\gamma_{j}\rangle\\
&=\sum_{i,j}\langle h,e_{i}\rangle\langle e_{i},k\gamma_{j}\rangle\langle \gamma_{j},h_{0}\rangle\\
&=\sum_{j}\langle h,k\gamma_{j}\rangle\langle \gamma_{j},h_{0}\rangle\\
&=\sum_{j}\langle k^{*}h,\gamma_{j}\rangle\langle \gamma_{j},h_{0}\rangle\\
&=\langle k^{*}h ,h_{0}\rangle
\end{align*} 
The proof is complete.
\end{proof}
 
 \begin{theorem}
 Let $ k \in \mathcal{B}(H_{0},H)$, and let  $f:X \to H$ be c-Bessel mappings for $H$, and $g:X\to H_0$ be a c-Bessel mapping for $H_0$.
 \item
If $k$ is onto then the following condition  is equivalent  with the assertions of the Theorem \ref{t3}:
\begin{align*}
\Vert kh_0\Vert^{2} = \int_{X}\langle h_0,g(x)\rangle \langle f(x),kh_0\rangle \,d\mu, \ \ h_0\in{H_0}.
\end{align*}  \item
 If $k^{*}$ is onto then the following condition  is equivalent  to the assertions of the Theorem \ref{t3}:

\begin{align*}
||k^{*}h|| = \int_{X}\langle h,f(x)\rangle \langle g(x),k^{*} h\rangle \,d\mu, \ \  h \in H.
\end{align*}
 \end{theorem}
 \begin{proof} (1)
Let $F:H_0 \to H_0$ be defined by
$$F(h_0):=T_g \langle kh_0,f \rangle, \ \ h_0 \in H_0.$$
$F$ is clearly linear and bounded, since for each $h_0 \in H_0,$

\begin{align*}
||F(h_0)||&= \sup_{k_0 \in {(H_o)}_{1}} |\langle F(h_0), k_0 \rangle |\\& =
\sup_{k_0 \in {(H_o)}_{1}} |\int_X  \langle g(x),k_0 \rangle \langle k h_0, f(x) \rangle d\mu|
 \\& \leq   \sup_{k_0 \in {(H_o)}_{1}} (\int_X |\langle k_0, g(x) \rangle | ^2 d \mu)^{1/2}  \sup_{\ell_0 \in {(H_o)}_{1}} (\int_X |\langle \ell_0 , f(x) \rangle |^2d \mu)^{1/2}||kh_0|| \\& \leq \sup_{k_0 \in {(H_o)}_{1}} (\int_X |\langle k_0, g(x) \rangle | ^2 d \mu)^{1/2}  \sup_{\ell_0 \in {(H_o)}_{1}} (\int_X |\langle \ell_0 , f(x) \rangle |^2d \mu)^{1/2}||k|| ||h_0||.
\end{align*}
For each $h_0 \in H_0$ we have
\begin{align*}
\langle  h_0, k^* kh_0 \rangle & = ||kh_0||^2 \\ & = \int_X \langle h_0, g(x) \rangle \langle f(x), k h_0 \rangle d \mu \\ & =\overline{ \langle T_g \langle  k h_0 , f \rangle , h_0 \rangle} \\& = \langle h_0 , T_g \langle kh_0, f \rangle \rangle.
\end{align*}
Hence, $ k^*kh_0= T_g \langle kh_0, f \rangle$. Since $k$ is onto, then we have  
$$ k^* h = T_g \langle h,f \rangle, \ \ h \in H.$$

The part (2) follows similarly.
 \end{proof}
\begin{definition}
 Let $ k \in \mathcal{B}(H_{0},H)$, and let  $f:X \to H$ be c-Bessel mappings for $H$, and $g:X\to H_0$ be a c-Bessel mapping for $H_0$.  We say that $f$,$g$  is a  c$k$-dual pair, if one of the assertions of the theorem \ref{t3} holds.
\end{definition}
\begin{lemma} Let $ k \in \mathcal{B}(H_{0},H)$, and let  $f:X \to H$ be c-Bessel mappings for $H$, and $g:X\to H_0$ be a c-Bessel mapping for $H_0$.
Let $f$, $g$ be ck-dual pair. Then $f$  is a ck-frame for $H$ with respect to $H_0$, and $g$ is a $ck^{*}$-frame for $H_0$ with respect to $H$.
\end{lemma}
\begin{proof}
 For each  $h\in H_0$ 
\begin{align*}
\Vert kh\Vert^4&=\vert \langle T_f(\langle h,g\rangle),kh\rangle\vert^2\\
&=\vert \int_{X}\langle h,g(x)\rangle.\langle f(x),kh\rangle\,d\mu\vert^2\\
&\leq(\int_{X}\vert \langle h,g(x)\rangle\vert^2\,d\mu)(\int_{X}\vert\langle f(x),kh\rangle\vert^2\,d\mu)\\
&\leq(\int_{X}\vert \langle h,g(x)\rangle\vert^2\,d\mu) B\Vert kh\Vert^2
\end{align*}
where $B$ is an upper bound of  the c-Bessel mapping $f$. This showes that $g$ is a c$k^{*}$-frame for $H_0$  with respect to $H$,  with lower bound $B^{-1}$. 

Similarly, $f$ is a c$k$-frame for $H$.
\end{proof}
In the following result, we construct some c$k$-duals with the same of Proposition 2.3 in \cite{arab}.
\begin{theorem}
Let $k\in\mathcal{B}(H_0,H)$  with closed range, and
 $f:X\rightarrow H$ be a c$k$-frame for $H$ with respect to $H_0$ and Bochner integrable (see \cite{rahmani}). Then $k^{\ast}(S_f|_{\mathcal{R}(k)})^{-1}\pi_{S_f(\mathcal{R}(k))}f$ is a ck-dual of $\pi_{\mathcal{R}(k)}f$  with the bounds $B^{-1}$ and $A^{-1}\Vert k\Vert^2 \Vert k^{\dagger}\Vert^2$, respectively, where $A$ and $B$ are ck-frame bounds for $f$.
\end{theorem}
\begin{proof}
By Theorem \ref{re1}, $S_f$ is invertible on $\mathcal{R}(k)$. 
Since $$k^{\ast}(S_f|_{\mathcal{R}(k)})^{-1}\pi_{S_f(\mathcal{R}(k))}\in \mathcal{B}(H,H_0),$$ then by Theorem \ref{t11},  $k^{\ast}(S_f|_{\mathcal{R}(k)})^{-1}\pi_{S_f(\mathcal{R}(k))}f$ is a c-Bessel mapping for $H_0$. Moreover, $(S_f|_{\mathcal{R}(k)})^{-1}S_{f}\vert_{\mathcal{R}(k)}=I_{\mathcal{R}(k)}$. Hence, for each $h\in H$
and $h_0\in H_0,$ we have
\begin{align*}
\langle h_0, k^{\ast}h \rangle&=\langle \big((S_f|_{\mathcal{R}(k)})^{-1}S_{f}\big)^{\ast}kh_0,h\rangle\\
&=\langle S^*_{f}\big((S_f|_{\mathcal{R}(k)})^{-1}\big)^{\ast}kh_0,h\rangle\\
&=\langle \pi_{\mathcal{R}(k)}S_f|_{S(\mathcal{R}(k))} \pi_{S_f(\mathcal{R}(k))}((S_f|_{\mathcal{R}(k)})^{-1})^{\ast}kh_0,h\rangle\\
&=\langle T_f\big(\langle \pi_{S_f(\mathcal{R}(k))}((S_f|_{\mathcal{R}(k)})^{-1})^{\ast} kh_0,f\rangle\big),\pi_{\mathcal{R}(k)}h\rangle\\
&=\int_{X}\langle \pi_{S_f(\mathcal{R}(k))}((S_f|_{\mathcal{R}(k)})^{-1})^{\ast}kh_0,f(x)\rangle\langle f(x),\pi_{\mathcal{R}(k)}h\rangle\,d\mu\\
&=\int_{X}\langle h_0,k^{\ast}(S_f|_{\mathcal{R}(k))}^{-1}\pi_{S_f(\mathcal{R}(k))}f(x)\rangle\langle \pi_{\mathcal{R}(k)}f(x),h\rangle d\mu.
\end{align*}
Therefore,
$$\langle  k^{\ast}h, h_0 \rangle=\int_{X} \langle h,\pi_{\mathcal{R}(k)}f(x)\rangle \langle k^{\ast}(S_f|_{\mathcal{R}(k)})^{-1}\pi_{S_f(\mathcal{R}(k))}f(x), h_0\rangle d\mu.$$
Now, put $g:=k^{\ast} (S_f|_{\mathcal{R}(k)})^{-1}\pi_{S_f(\mathcal{R}(k))}f;$ by Theorem   \ref{t3}, $g$ is a c$k$-dual for $\pi_{\mathcal{R}(k)}f$ with lower bound $B^{-1}$. Furthermore, by Theorem  \ref{re1}, for each $h\in\mathcal{R}(k),$ 
\begin{align*}
\Vert ((S_f|_{\mathcal{R}(k)})^{-1})^{\ast}h\Vert^2&=\langle S_f|_{\mathcal{R}(k)}^{-1}(S_f|_{\mathcal{R}(k)}^{-1})^{\ast}h ,h\rangle\\
&\leq A^{-1}\Vert k^{\dagger}\Vert^2 \Vert ((S_f|_{\mathcal{R}(k)})^{-1})^{\ast}h\Vert\Vert h\Vert.
\end{align*} 
Therefore, for each $h\in H$
\begin{align*}
\int_{X}\vert\langle h,g(x)\rangle\vert^2\,d\mu&=\int_{X}\vert\langle \pi_{S_f(\mathcal{R}(k))}\big((S_f|_{\mathcal{R}(k)}\big)^{-1})^{\ast}kh,f(x)\rangle\vert^2\,d\mu\\
&\leq\langle S_f|_{\mathcal{R}(k)}((S_f|_{\mathcal{R}(k)})^{-1})^{\ast}kh,((S_f|_{\mathcal{R}(k)})^{-1})^{\ast}kh\rangle\\
&=\langle \big((S_f|_{\mathcal{R}(k)})^{-1}\big)^{\ast}kh, kh\rangle\\
&\leq\Vert kh\Vert\Vert \big((S_f|_{\mathcal{R}(k)})^{-1}\big)^{\ast}kh\Vert\\
&\leq\Vert kh\Vert A^{-1}\Vert k^{\dagger}\Vert^2 \Vert kh\Vert\\
&\leq A^{-1}\Vert k\Vert^2\Vert k^{\dagger}\Vert^2\Vert h\Vert^2
\end{align*}
and the result follows.
\end{proof}



\begin{thebibliography}{20}

\bibitem{arab} F. Arabyani Neyshaburi,  and A. A. Arefijamaal,
Some constructions of $K$-frames and their duals.
\emph{Rocky Mountain Journal of Math}. \textbf{47}(6), (2017), 1749-1764.

\bibitem {am00} J. Benedetto, A. Powell and O. Yilmaz,
 Sigma-Delta quantization and finite frames,
\emph{ IEEE Trans. Inform. Th.} 52, (2006) 1990-2005.

\bibitem {aklr}  H. Bolcskel , F. Hlawatsch and H.G Feichyinger,
 Frame-Theoretic analysis of oversampled filter bank,
\emph{ IEEE Trans. Signal Processing.} \textbf{46}(12), (1998) 3256- 3268.

\bibitem {bbk}  E. J. Candes and D. L. Donoho,
 New tight frames of curvelets and optimal representation of
 objects with piecwise $C^2$ singularities,
\emph{ Comm. Pure and App. Math.} 56 (2004), 216-266.

\bibitem {cay}  P. G. Casazza and   G. Kutyniok,
 Frame of subspaces,
\emph{ Wavelets, Frames and Operator Theory.} (College Park, MD, 2003), Contemp. Math. 345, Amer. Math. Soc., Providence, RI, (2004) , 87-113.

\bibitem {cz02} P. G. Casazza ,   G. Kutyniok, and  S. Li,
 Fusion frames and Distributed Processing,
\emph{ Appl. Comput. Harmon. Anal.} 25 (2008), 114-132.

\bibitem {dt97}  P. G. Casazza and J,  Kova\^ cevi\' c,
 Equal-norm tight frames with erasures,
\emph{ Adv. Comput. Math.} 18 (2003), 387-430.

\bibitem {gj91} O. Christensen,
\emph{ Introduction to frames and Riesz bases,}
Boston, Birkhauser.  2003.

\bibitem{ga94} I. Daubechies, A. Grossmann and Y. Meyer, 
 Painless nonorthogonal Expansions,
\emph{J. Math. Phys.} 27 (1986), 1271-1283.

\bibitem {do} R. G. Douglas,
On majorization, factorization and range inclusion of operators on Hilbert spaces,
\emph{ Proc. Amer. Math. Soc.} \textbf{17}(2) (1966), 413-415.

\bibitem {hy44} R. J. Duffin and  A. C. Schaeffer, 
 A class of nonharmonik Fourier series,
\emph{ Trans. Amer. Math. Soc.} 72 (1952), 341-366.

\bibitem {hy41} M. H. Faroughi and  E. Osgooei, 
 C-Frames and C-Bessel Mappings,
\emph{ Bulletin of IMS,} \textbf{38}(1) (2012), 203-222.

\bibitem {hy43} J. P. Gabardo and D. Han,
 Frames associated with measurable space,
\emph{ Adv. comp. Math.} 18 (2003), 127-147.

\bibitem {hy143} L. G\u avru\c ta, L,  
 Frames for operators,
\emph{ Appl. Comput. Harmon. Anal.} 32 (2012), 139-144. 

\bibitem {hir98} B.  Hassibi, B.  Hochwald,  A. Shokrollahi and  W. Sweldens,
Representation theory for high-rate multiple-antenna code design,
\emph{IEEE Trans. Inform. Theory.} 47 (2001), 2335-2367.

\bibitem {hir103} A. Rahimi A. Najati and Y. N. Dehgan,
Continuous frame in Hilbert space,
\emph{ Methods of Functional Analysis and Topology,} 12 (2006), 170-182.

\bibitem {hir1103} A. Rahimi, A. Najati and M. H. Faroughi,
 Continuous and discrete frames of subspaces in Hilbert spaces,
\emph{ Southeast Asian Bulletin of Mathematics,} 32 (2008), 305-324.

\bibitem {rahmani} M. Rahmani,
On some properties of c-frames,
\newblock J. Math. Research with Appl., \textbf{37}(4) (2017), 466-476.

\bibitem {ir96} W. Rudin,
\emph{ Functional Analysis,}
 New York, Tata Mc Graw-Hill Editions, 1973.

 \bibitem {ir97} W. Rudin, 
\emph{ Real and Complex Analysis,}
 \newblock New York, Tata Mc Graw-Hill Editions, 1987.

\bibitem {ir98} S. Sakai, 
\emph{ $C^*$-Algebras and $W^*$-Algebras,}
 New York, Springer-Verlag, 1998.

\bibitem {XX} X. Xiao, Y. Zhu and L. Gavruta, 
 Some Properties of $K$-frames in Hilbert Spaces,
\emph{ Results. Math.} 63,  (2012), 1243-1255.
\end{thebibliography}
\end{document}